\documentclass[journal,a4paper]{IEEEtran}
\usepackage[T1]{fontenc}% optional T1 font encoding
\usepackage{layouts}
\usepackage{cite}
\usepackage[pdftex,colorlinks=true,bookmarks=true,citecolor=blue,urlcolor=blue]{hyperref} %pdflatex

\ifCLASSINFOpdf
  \usepackage[pdftex]{graphicx}
   \graphicspath{{./figs/}}
  \DeclareGraphicsExtensions{.pdf,.jpeg,.png}
\else
\fi
\usepackage{amsmath}
\usepackage{amsthm}
\usepackage[cmintegrals]{newtxmath}
\usepackage{bm}
\usepackage{algorithmic}
\usepackage{array}
\ifCLASSOPTIONcompsoc
  \usepackage[caption=false,font=normalsize,labelfont=sf,textfont=sf]{subfig}
\else
  \usepackage[caption=false,font=footnotesize]{subfig}
\fi
\usepackage{mathrsfs}  
\usepackage{amsbsy}

\newcommand{\field}[1]{\mathbb{#1}}
\newcommand{\fs}[1]{\mathsf{#1}}

\DeclareMathOperator{\diag}{diag}

\newcommand{\tp}{\intercal}% transpose operation

\newcommand{\bigO}[1]{\mathop{\mathcal{O}}\left(#1\right)}

\let\Re\relax
\DeclareMathOperator{\Re}{Re}

\newcommand{\vv}[1]{\mathbf{#1}}
\newcommand{\vs}[1]{\boldsymbol{#1}}

\DeclareMathOperator{\sech}{sech}

\newcommand{\wtilde}[1]{\widetilde{#1}}

\newcommand{\et}{\textit{et~al.}}

\usepackage{enumitem}

\begin{document}
\title{Efficient Nonlinear Fourier Transform Algorithms of Order Four on Equispaced Grid}

%\author{Michael~Shell,~\IEEEmembership{Member,~IEEE,}
%        John~Doe,~\IEEEmembership{Fellow,~OSA,}
%        and~Jane~Doe,~\IEEEmembership{Life~Fellow,~IEEE}% <-this % stops a space

\author{Vishal Vaibhav% <-this % stops a space
\thanks{Email:~\tt{vishal.vaibhav@gmail.com}}% <-this % stops a space
}

\IEEEpubid{\begin{minipage}[t]{\textwidth}
\vskip2em
\centering
\copyright~2019 IEEE. Personal use of this material is permitted. Permission 
from IEEE must be obtained for all other uses, including reprinting/republishing this 
material for advertising or promotional purposes, collecting new collected works for 
resale or redistribution to servers or lists, or reuse of any copyrighted component 
of this work in other works.
\end{minipage}}


\maketitle

\begin{abstract}
We explore two classes of exponential integrators in this letter to design
nonlinear Fourier transform (NFT) algorithms with a desired
accuracy-complexity trade-off and a convergence order of $4$ on an 
equispaced grid. The integrating factor based method in the class of
Runge-Kutta methods yield algorithms with complexity $O(N\log^2N)$ (where $N$ is
the number of samples of the signal) which have superior accuracy-complexity
trade-off than any of the fast methods known currently.
The integrators based on Magnus series expansion, namely, standard and 
commutator-free Magnus methods yield algorithms of complexity $O(N^2)$ that 
have superior error behavior even for moderately small step-sizes and higher 
signal strengths.
\end{abstract}

%\begin{IEEEkeywords}
%\end{IEEEkeywords}

% \fi
\IEEEpeerreviewmaketitle

%%%%%%%%%%%%%%%%%%%%%%%%%%%%%%%%%%%%%%%%%%%%%%%%%%%%%%%%%%%%%%%%%%%%%%%%%%%%%%%%%%%
\section{Introduction}
In a series of papers~\cite{V2017INFT1,V2018BL,V2018LPT}, it
was shown recently that exponential linear multistep methods (LMMs) provide 
a natural setting for higher-order convergent fast nonlinear Fourier 
transform (NFT). This followed from a simple observation
that the transfer matrices obtained are amenable 
to FFT-based fast polynomial arithmetic. Note that in the earliest works on fast 
NFTs~\cite{WP2013d}, the Ablowitz-Ladik method can be interpreted as the 
exponential Euler method. In this paper, we use the exponential 
Runge-Kutta methods to obtain a family of fast NFTs (provided 
that the nodes are equispaced). In particular, we
present two fast NFTs with fourth order accuracy based on fourth order
Runge-Kutta methods. The structure of the transfer matrix reveals that such methods are
superior to those based on LMMs in terms of complexity while the numerical tests
reveal that they also have a superior accuracy-complexity trade-off. The first 
algorithm is based on the classical 
fourth order (explicit) Runge-Kutta method which has
been studied by several authors in the context of 
NFTs~\cite{BCT1998,CPW2018TSP}. The second method uses the three-stage 
Lobatto~IIIA (implicit) Runge-Kutta method.

\begin{figure*}[!ht]
\centering
\includegraphics[scale=1]{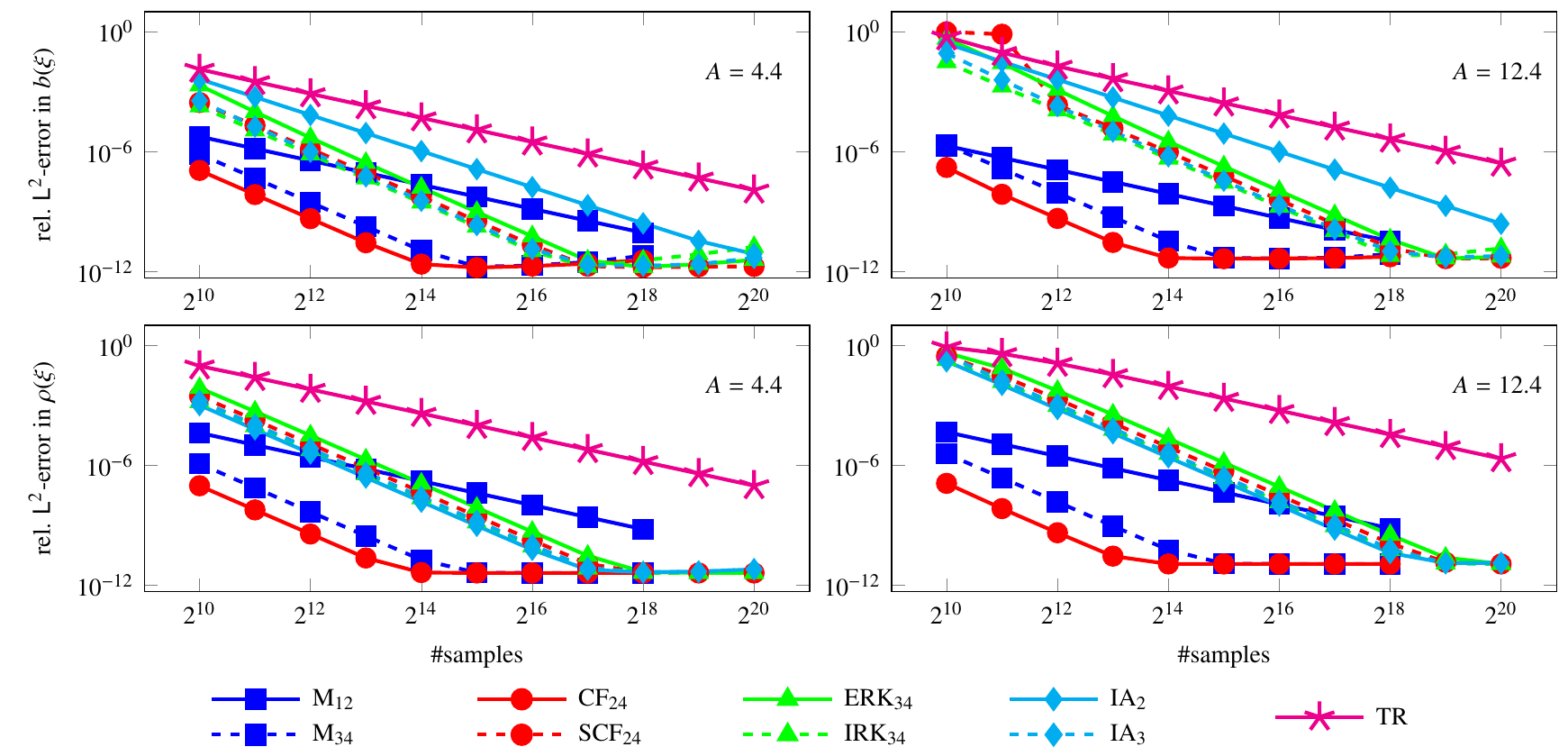}%
\caption{\label{fig:convg}The figure shows the convergence analysis of 
various methods for the secant-hyperbolic profile. It is evident that the CF
method with two exponentials (CF$_{24}$) outperforms fourth order Magnus method
(M$_{34}$). Note that the complexity of methods M$_{12}$, M$_{34}$ and CF$_{24}$
is $O(N^2)$ while the rest are of $O(N\log^2N)$.}
%\end{figure*}
%\begin{figure*}[!ht]
\vskip0.5em
%\centering
\includegraphics[scale=1]{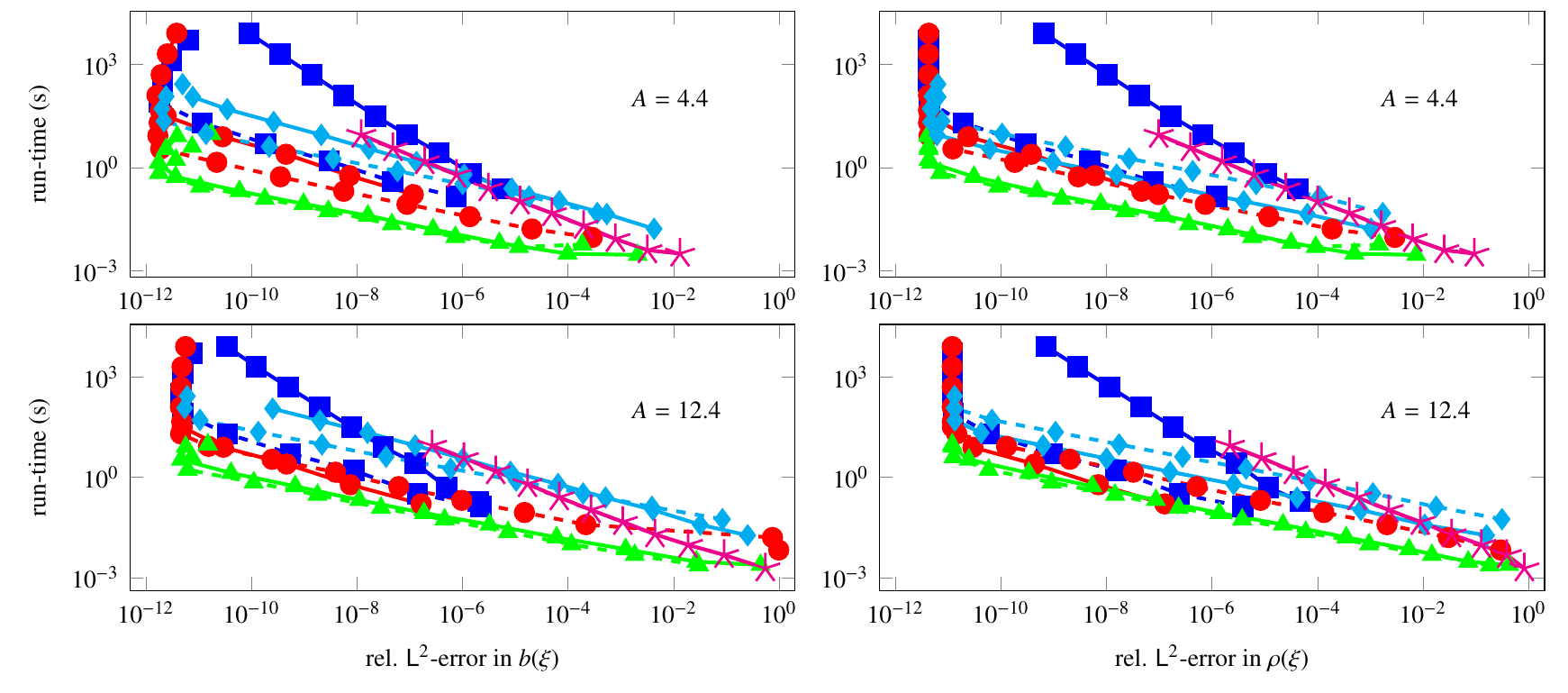}%
\caption{\label{fig:trade-off}The figure shows the 
accuracy-complexity trade-off for various methods
for the secant-hyperbolic profile. The legends are same
as that of Fig.~\ref{fig:convg}. The plots demonstrate that the Runge-Kutta
methods (ERK$_{34}$, IRK$_{34}$) turn out to be far superior than all the other
`fast' methods with regard to the accuracy-complexity trade-off. It is
interesting to note that the `slow' methods become quite competitive with
increasing signal strength.}
\end{figure*}

For moderately small step-sizes, most fast methods yield poor accuracy 
specially corresponding to the large values of the spectral parameter. An error 
analysis of such integrating factor based methods~\cite{BCT1998} shows that 
the error terms contain positive powers of the spectral parameter which
necessitates the use of smaller step-sizes in order to keep the error low. On the
contrary, integrators based on Magnus series expansion are known to have 
error terms that contain negative powers of the spectral parameter. In this
letter, we follow the recipe presented by Blanes~\et~\cite{BCR2000} to develop 
a fourth order Magnus method and a fourth order commutator-free (CF) 
Magnus method~\cite{BM2006}, both of which take samples of the potential on an
equispaced grid in order to compute the NFT. Let us emphasize that our CF 
method is different from those
considered in~\cite{CPW2018} where the samples of the potential are needed on
the Gauss-Legendre nodes (the authors generate the samples by interpolation on 
an equispaced grid, locally). The fourth-order (standard) Magnus method happens to be 
faster than the corresponding CF method on account of the fact that there is an 
additional matrix exponential introduced in the CF method in order to avoid the 
use of commutators. The accuracy-complexity trade-off, however, is similar 
for the two methods. Finally, we also present a fast variant of the CF 
method of order four (formally) by employing the fourth order splitting on the lines 
of~\cite{D2001,CPW2018TSP}. Despite the well-known limitation imposed on the order
and stability of such techniques as demonstrated by Sheng~\cite{S1989}, we do
not find any reduction of order within the double precision arithmetic. It is
noteworthy that the authors in~\cite{CPW2018TSP} found the aforementioned 
splitting worsen in accuracy after a certain step-size. It is also not clear
from their analysis if this scheme is convergent in their setting.

We begin our discussion with a brief review of the scattering theory closely
following the formalism presented in~\cite{AKNS1974}. The nonlinear Fourier
transform of any signal is defined via the Zakharov-Shabat (ZS) scattering 
problem which can be stated as follows: For $\zeta\in\field{R}$ and $\vv{v}=(v_1,v_2)^{\tp}$,
\begin{equation}\label{eq:ZS-prob}
\vv{v}_t =\left[-i\zeta\sigma_3+U(t)\right]\vv{v}
\equiv\mathcal{T}(t;\zeta)\vv{v},
\end{equation}
where $\sigma_3=\diag(1,-1)$. The potential 
$U(t)$ is defined by $U_{11}=U_{22}=0,\,U_{12}=q(t)$ and $U_{21}=r(t)$ with $r=\kappa q^*$
($\kappa\in\{+1, -1\}$). Here, $\zeta\in\field{R}$ is known as the \emph{spectral parameter}
and $q(t)$ is the complex-valued signal. The solution of the scattering 
problem~\eqref{eq:ZS-prob}, henceforth referred
to as the ZS problem, consists in finding the so called 
\emph{scattering coefficients} which are defined through special solutions 
of~\eqref{eq:ZS-prob} known as the \emph{Jost solutions} which are linearly 
independent solutions of~\eqref{eq:ZS-prob} such that they have a plane-wave 
like behavior at $+\infty$ or $-\infty$. The Jost solution of the 
\emph{second kind}, denoted by $\vs{\phi}(t,\zeta)$, has 
the asymptotic behavior $\vs{\phi}(t;\zeta)e^{i\zeta t}\rightarrow(1,0)^{\tp}$ 
as $t\rightarrow-\infty$. The asymptotic behavior 
$\vs{\phi}(t;\zeta)\rightarrow (a(\zeta)e^{-i\zeta t}, b(\zeta)e^{i\zeta
t})^{\tp}$ as $t\rightarrow\infty$ determines the scattering coefficients $a(\zeta)$ and
$b(\zeta)$ for $\zeta\in\field{R}$. In this letter, we primarily focus on 
the continuous spectrum, also referred to as the 
\emph{reflection coefficient}, which is defined by 
$\rho(\xi)={b(\xi)}/{a(\xi)}$ for $\xi\in\field{R}$.

\section{The Numerical Scheme}\label{eq:num-scheme}
\subsection{Runge-Kutta Method}
In this section, we will develop the integrating factor based exponential
Runge-Kutta (RK) method for the numerical solution of the ZS 
problem. Following~\cite{V2017INFT1,V2018LPT}, 
we begin with the transformation $\tilde{\vv{v}}=e^{i\sigma_3\zeta t}\vv{v}$
so that~\eqref{eq:ZS-prob} becomes $\tilde{\vv{v}}_t=\wtilde{U}\tilde{\vv{v}}$
with $\wtilde{U}=e^{i\sigma_3\zeta t}Ue^{-i\sigma_3\zeta t}$ whose entries are 
$\wtilde{U}_{11}=\wtilde{U}_{22}=0,\,\wtilde{U}_{12}=q(t)e^{2i\zeta t}$ and
$\wtilde{U}_{21}=r(t)e^{-2i\zeta t}$. 
Let the step size be $h>0$ and the quantities $c_j\in[0,1]$ be ordered 
so that the nodes within the step can be stated as 
$t_n\leq t_n+c_1h\leq t_n+c_2h\leq\ldots\leq t_n+c_sh\leq t_{n+1}$. For the potential
sampled at these nodes, we use the convention
$Q_{n+c_k}=hq(t_n+c_kh)$, $R_{n+c_k}=hr(t_n+c_kh)$ and
$\wtilde{U}_{n+c_k}=\wtilde{U}(t_n+c_kh)$.
In order for the resulting discrete system to be amenable to FFT-based fast polynomial
arithmetic, it is sufficient to have each of the $c_i$'s belong to the 
set of uniformly distributed nodes in $[0,1]$. A $s$-stage RK method is 
characterized by the nodes
$\vv{c}=(c_1,c_2,\ldots,c_s)\in\field{R}^s$, and, the weights 
$\vv{b}=(b_1,b_2,\ldots,b_s)\in\field{R}^s$ and
$(a_{ij})\in\field{R}^{s\times s}$. Introducing the intermediate stage quantities
$\tilde{\vv{v}}_{n,k}$ for $k=1,2,\ldots,s$, we have
\begin{equation}
\left\{\begin{aligned}
&\tilde{\vv{v}}_{n,j}=\tilde{\vv{v}}_{n}
+h\sum_{k=1}^sa_{jk}\wtilde{U}_{n+c_k}\tilde{\vv{v}}_{n,k},\quad
j=1,2,\ldots,s,\\
&\tilde{\vv{v}}_{n+1}=\tilde{\vv{v}}_{n}
+h\sum_{k=1}^sb_{k}\wtilde{U}_{n+c_k}\tilde{\vv{v}}_{n,k}.
\end{aligned}\right.
\end{equation}
This system of equations can be solved in any computer algebra system to obtained
the transfer matrix connecting the vectors $\vv{v}_{n+1}$ to $\vv{v}_n$.

Setting $z=\exp(i\zeta h/2)$, the 
Lobatto~IIIA method (labelled as IRK$_{34}$) of order $4$~\cite{Butcher2003} simplifies to
\begin{align}
{\vv{v}}_{n+1}&={z^{-2}}\left[\Delta_{n+1}(z^2)\right]^{-1}M_{n+1}(z^2){\vv{v}}_{n},
\label{eq:LobattoIIIA4}\\
M_{n+1}
&=\begin{pmatrix}
1+\frac{z^2}{12}Q_{n+1}R_{n+1/2} &\frac{z^2}{6}Q_{n+1}+\frac{1}{3}Q_{n+1/2}\\
\frac{1}{6}R_{n+1}+\frac{z^2}{3}R_{n+1/2}& z^2+\frac{1}{12}R_{n+1}Q_{n+1/2}
\end{pmatrix}\nonumber\\
&\quad\times\begin{pmatrix}
1+\frac{z^2}{12}R_nQ_{n+1/2}&\frac{1}{6}Q_{n}+\frac{z^2}{3}Q_{n+1/2}\\
\frac{z^2}{6}R_{n}+\frac{1}{3}R_{n+1/2}&z^2+\frac{1}{12}Q_nR_{n+1/2}
\end{pmatrix},\label{eq:LIIIA3-TM}\\
\Delta_{n+1}
&=\left(1+\frac{z^{-2}}{12}R_{n+1}Q_{n+1/2}\right)
\left(1+\frac{z^2}{12}Q_{n+1}R_{n+1/2}\right)\nonumber\\
&\quad-\frac{1}{36}\left(Q_{n+1}+2{z^{-2}}Q_{n+1/2}\right)
\left(R_{n+1}+2{z^2}R_{n+1/2}\right).
\end{align}
The fourth order classical RK method (labelled as ERK$_{34}$)
simplifies to the form~\eqref{eq:LobattoIIIA4} with $\Delta_{n+1}=1$
and the entries of the transfer matrix are given by
\begin{equation}
\left\{\begin{aligned}
M^{(n+1)}_{11}(z^2)&= 
G_{n+1/2}+\frac{z^2}{6}\left(Q_{n+1/2}R_{n}+Q_{n+1}R_{n+1/2}\right)\\
&\quad+\frac{z^4}{24}Q_{n+1/2}Q_{n+1}R_nR_{n+1/2},\\
M^{(n+1)}_{22}(z^2)&=\frac{1}{24}Q_{n}Q_{n+1/2}R_{n+1/2}R_{n+1}+G_{n+1/2}z^4\\
&\quad+\frac{z^2}{6}\left(Q_{n}R_{n+1/2} + Q_{n+1/2}R_{n+1}\right),\\
M^{(n+1)}_{12}(z^2)
&=\left(\frac{z^2}{6}Q_{n}+\frac{z^4}{6}Q_{n+1}\right)H_{n+1/2}
+\frac{2}{3}Q_{n+1/2},\\
M^{(n+1)}_{21}(z^2)&=
\left(\frac{z^4}{6}R_{n}+\frac{1}{6}R_{n+1}\right)H_{n+1/2}
+\frac{2z^2}{3}R_{n+1/2},
\end{aligned}\right.
\end{equation}
where $6G_{n+1/2}-6=Q_{n+1/2}R_{n+1/2}=2H_{n+1/2}-2$.

\begin{figure}[!ht]
\centering
\includegraphics[scale=1]{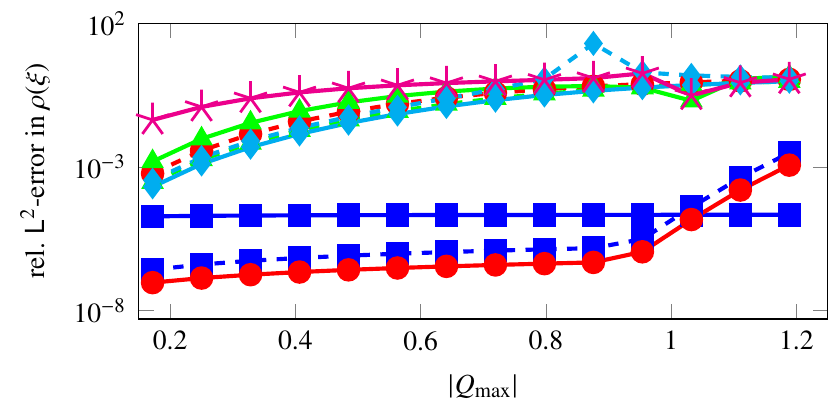}%
\caption{\label{fig:err}The figure shows the error behavior of the 
NFT algorithms for $N=2^{10}$
as a function of $Q_{\text{max}}=h\|q\|_{\infty}$ for the secant-hyperbolic
profile. The legends are same as that
of Fig.~\ref{fig:convg}. It is straightforward to conclude that for smaller number of 
samples, the `slower' methods (M$_{12}$, M$_{34}$ and CF$_{24}$) are far
superior to the `fast' methods for $\xi\in\Omega_h$, the
principal domain of a given method, with comparable run-times.}
\end{figure}

\subsubsection{Scattering coefficients}
Let the computational domain be $\Omega=[T_1,T_2]$ and set $2T=T_2-T_1$. Let the number
of steps be $N_s$ and $h=2T/N_s$. Also, let $h\ell_+=T_2$ and 
$h\ell_-=-T_1$. The grid is defined by $t_n= T_1 + nh/2,\,\,n=0,1,\ldots,N,$ with 
$t_{N}=T_2$ where $N=2N_s$ is the number of samples. Let the 
potential $q(t)$ be supported in $\Omega$ and we assume $q_0=0$ for convenience. In 
order to represent the Jost solutions, we introduce the polynomial vector
\begin{equation}\label{eq:poly-vec-rk}
\vv{P}_n(z)
=\begin{pmatrix}
P^{(n)}_{1}(z)\\
P^{(n)}_{2}(z)
\end{pmatrix}
=\sum_{j=0}^{2n}
\vv{P}^{(n)}_{j}z^j
=\sum_{j=0}^{2n}
\begin{pmatrix}
P^{(n)}_{1,j}\\
P^{(n)}_{2,j}
\end{pmatrix}z^{j},
\end{equation}
and the polynomial $D_n(z) = \sum_{j=0}^{2n}D^{(n)}_jz^j$. Consider 
the Jost solution $\phi(t;\zeta)$. For the Lobatto IIIA method, the 
Jost solution can be stated as 
$\vs{\phi}_n(z^2) = (z^2)^{\ell_-}\left[D_{n}(z^2)\right]^{-1}{\vv{P}_{n}(z^2)}$,
with
\begin{equation}
\vv{P}_{n+1}(z^2)=\Theta^{-1}_{n+1}M_{n+1}(z^2)\vv{P}_{n}(z^2),
\end{equation}
where $\Theta_{n}$ is the constant part of $\Delta_n(z^2)$ and
\begin{equation}
D_{n+1}(z^2)=\Theta^{-1}_{n+1}z^2\Delta_{n+1}(z^2)D_{n}(z^2).
\end{equation}
From 
${\vs{\phi}}(T_2;\zeta)=\left(ae^{-i\zeta T_2},
be^{+i\zeta T_2}\right)^{\tp}$, it follows that the discrete scattering
coefficients are given by
\begin{equation}
\begin{split}
&a_{N_s}(z^2) = (z^2)^{\ell_-+\ell_+}{D^{-1}_{N_s}(z^2)}{P^{(N_s)}_1(z^2)},\\
&b_{N_s}(z^2) = (z^2)^{\ell_--\ell_+}{D^{-1}_{N_s}(z^2)}{P_2^{(N_s)}(z^2)},
\end{split}
\end{equation}
where $z=\exp(i\zeta h/2)$. For the classical fourth order RK method, 
$\Delta_n(z^2)=\Theta_n=1$; therefore, $a_{N_s}(z^2) = {P^{(N_s)}_1(z^2)}$ and 
$b_{N_s}(z^2) = (z^2)^{-2\ell_+}{P_2^{(N_s)}(z^2)}$.
The principal branch for the discrete scattering coefficients here works out to
be $\Re\zeta\in[-\pi/2h,\pi/2h]$. This follows from the principle branch of the
individual transfer matrices. The nodes
$\zeta_j=\xi_j+i0=j(\pi/h)/N,\,j\in\field{Z},$ lead to
$z_j^2=\exp(i\pi j/N)$ which is not in the standard form for FFT algorithms to be
used. Therefore, we would like to work with $N'=2N$ nodes so that 
$z_j^2=\exp(i2\pi j/N')$ and pad the input vector with zeros. Following as 
in~\cite{V2018LPT}, the complexity of
computing the scattering coefficient works out to be $O(N\log^2N)$.

\subsection{Standard and commutator-free Magnus methods}
\subsubsection{Magnus method} 
Let us assume that the solution of the ODE~\eqref{eq:ZS-prob} can be
written as $\vv{v}(t)=\exp[\Lambda(t;t_n)]\vv{v}(t_n)$ for 
$t\in[t_{n+1},t_n]$, then $\Lambda(t;t_n)$ has a series representation 
known as the \emph{Magnus series}~\cite{BCR2000}. Truncating this series to 
achieve the desired order of accuracy yields a family
of numerical schemes known as Magnus method. For the ZS problem, Magnus 
integrators preserve the Lie group structure of the Jost solution and its accuracy 
does not worsen with increasing $|\zeta|$. In designing a fourth order Magnus 
integrator (labelled as M$_{34}$), we follow the method due to
Blanes~\et~\cite{BCR2000,BM2006}: Defining
\begin{equation}
\mathcal{T}_{n+1}^{(j)}(h;\zeta)=h\int_{0}^{1}\left(\tau-\frac{1}{2}\right)^j
\mathcal{T}(t_{n}+\tau h;\zeta)d\tau,
\end{equation}
the method proceeds by expanding $\Lambda(t_{n+1};t_n)$ in terms of the quantities
$\mathcal{T}_{n+1}^{(j)}(h;\zeta)$ using the Magnus series. For the fourth order
method, setting $\Lambda(t_{n+1};t_n)\approx\Lambda_{n+1}$, we have~\cite{BCR2000,BM2006}
$\Lambda_{n+1}=\mathcal{T}_{n+1}^{(0)}
+\left[\mathcal{T}_{n+1}^{(1)},\mathcal{T}_{n+1}^{(0)}\right]$. Evaluating 
$\mathcal{T}_{n+1}^{(0)}$ and $\mathcal{T}_{n+1}^{(1)}$
upto fourth order accuracy using the three-point Gauss quadrature involving 
Legendre-Gauss-Lobatto (LGL) nodes, $\vv{c}=(0,{1}/{2}, 1)$, the numerical scheme
for the ZS problem can be stated as 
$\vv{v}_{n+1}=\exp(\Lambda_{n+1})\vv{v}_n$ where, 
% \begin{equation*}
% \Lambda_{n+1}=h\left(-i\zeta\sigma_3+\mathcal{U}_{n+1}\right)
% +\frac{h^2}{12}[U_{n+1}-U_{n},-i\zeta\sigma_3+\mathcal{U}_{n+1}],
% \end{equation*}
% where $\mathcal{U}_{n+1}=\frac{1}{6}(U_{n}+4U_{n+1/2}+U_{n+1})$. 
\begin{multline}
\Lambda_{n+1}=
\begin{pmatrix}
\Xi_{n+1} & G_{n+1}\\
H_{n+1} &-\Xi_{n+1}
\end{pmatrix}
+i\zeta h
\begin{pmatrix}
-1 & \frac{(Q_{n+1}-Q_{n})}{6}\\
\frac{(R_{n}-R_{n+1})}{6} & 1
\end{pmatrix},
\end{multline}
where $6G_{n+1}=(Q_{n}+4Q_{n+1/2}+Q_{n+1})$, 
$H_{n+1}=\kappa G^*_{n+1}$ and 
$12\Xi_{n+1}=[(Q_{n+1}-Q_{n})H_{n+1}-(R_{n+1}-R_{n})G_{n+1}]$.
%From here, it is easy to conclude that $\Lambda_{n+1}$ is traceless so that 
%$\det(e^{\Lambda_{n+1}})=1$. 
%Putting $\Gamma_{n+1} =
%\pm\sqrt{-\det(\Lambda_{n+1})}$, we have
%\begin{equation}\label{eq:expm}
%e^{\Lambda_{n+1}}=\cosh\left(\Gamma_{n+1}\right)\sigma_0
%+\left[{\sinh\left(\Gamma_{n+1}\right)}/{\Gamma_{n+1}}\right]\Lambda_{n+1}.
%\end{equation}
For the purpose of comparison, we would also like to consider the Magnus
method with one point Gauss quadrature (labelled as M$_{12}$) which 
is of order $2$~\cite{V2017INFT1}.

\subsubsection{Commutator-free Magnus method}
Blanes and Moan~\cite{BM2006} have constructed fourth-order commutator-free
(CF) methods that are based on Magnus method. Using the quantities defined above, 
the CF method (labelled as CF$_{24}$) can be stated as 
\begin{equation}
\vv{v}_{n+1}=e^{\Lambda^{(+)}_{n+1}}e^{\Lambda^{(-)}_{n+1}}\vv{v}_n,
\quad 2\Lambda^{(\pm)}_{n+1}=\mathcal{T}_{n+1}^{(0)}\pm 4\mathcal{T}_{n+1}^{(1)}.
\end{equation}
Evaluating $\mathcal{T}^{(0)}_{n+1}$ and $\mathcal{T}^{(1)}_{n+1}$
upto fourth-order accuracy using the three-point Gauss quadrature involving LGL
nodes, we have $2\Lambda^{(\pm)}_{n+1}=-{i\zeta h\sigma_3}+2h\mathcal{U}^{(\pm)}_{n+1}$
where $[\mathcal{U}^{(\pm)}_{n+1}]_{11}=[\mathcal{U}^{(\pm)}_{n+1}]_{22}=0$,
\begin{equation}
\begin{split}
&h[\mathcal{U}^{(+)}_{n+1}]_{12}\equiv
G^{(+)}_{n+1}=(3Q_{n+1}+4Q_{n+1/2}-Q_{n})/12,\\
&h[\mathcal{U}^{(-)}_{n+1}]_{12}\equiv
G^{(-)}_{n+1}=(3Q_{n}+4Q_{n+1/2}-Q_{n+1})/12,\\
\end{split}
\end{equation}
and $h[\mathcal{U}^{(\pm)}_{n+1}]_{21}
\equiv H^{(\pm)}_{n+1}=\kappa h[\mathcal{U}^{(\pm)}_{n+1}]^*_{12}$. Given that 
there are two matrix exponentials involved, the CF Magnus method has higher complexity than
that of the standard Magnus method. 

\subsubsection{A fast variant} 
The CF Magnus method further allows us to obtain a fast NFT algorithm (labelled
as SCF$_{24}$) via splitting of the matrix exponential. Consider a formally fourth order 
splitting~\cite{D2001}:
\begin{multline}
3e^{\Lambda^{(\pm)}_{n+1}}=4\left(e^{-\frac{1}{8}i\zeta
h\sigma_3}e^{\frac{1}{2}h\mathcal{U}^{(\pm)}_{n+1}}
e^{-\frac{1}{8}i\zeta h\sigma_3}\right)^2\\
-e^{-\frac{1}{4}i\zeta h\sigma_3}e^{h\mathcal{U}^{(\pm)}_{n+1}}
e^{-\frac{1}{4}i\zeta h\sigma_3}
+\bigO{h^5}.
\end{multline}
This splitting is stable and convergent~\cite{D2001}; however, 
its global order of convergence is $\leq2$~\cite{S1989}. 
Introducing
$16\Delta^{(\pm)}_{n+1}=12-3G_{n+1}^{(\pm)}H_{n+1}^{(\pm)}$ and 
$\Theta^{(\pm)}_{n+1}=1-G^{(\pm)}_{n+1}H^{(\pm)}_{n+1}$, and, 
putting $z=\exp(i\zeta h/4)$, the transfer matrix relation can be written as
\begin{equation}
\vv{v}_{n+1} ={z^{-4}}
\left[\Delta^{(+)}_{n+1}\Delta^{(-)}_{n+1}\right]^{-1}
M^{(+)}_{n+1}(z)M^{(-)}_{n+1}(z)\vv{v}_n,
\end{equation}
where the entries of the matrix $M^{(\pm)}_{n+1}(z)$ are
\begin{equation}
\left\{\begin{aligned}
&M^{(\pm,n+1)}_{11}(z)= \left(1 - C^{(\pm)}_{n+1}\right)
+\frac{z^2}{4}G_{n+1}^{(\pm)}H_{n+1}^{(\pm)},\\
&M^{(\pm,n+1)}_{22}(z)= \left(1-C^{(\pm)}_{n+1}\right)z^4 
+\frac{z^2}{4}G_{n+1}^{(\pm)}H_{n+1}^{(\pm)},\\
&M^{(\pm,n+1)}_{12}(z)=\frac{z}{2}\left(1+z^2\right)G_{n+1}^{(\pm)}
-z^2C^{(\pm)}_{n+1}G_{n+1}^{(\pm)},\\
&M^{(\pm,n+1)}_{21}(z)=\frac{z}{2}\left(1+z^2\right)H_{n+1}^{(\pm)}
-z^2C^{(\pm)}_{n+1}H_{n+1}^{(\pm)},
\end{aligned}\right.
\end{equation}
where
$3C^{(\pm)}_{n+1}={\Delta^{(\pm)}_{n+1}}[\Theta^{(\pm)}_{n+1}]^{-1/2}$. The 
discrete scattering coefficients can be written as $a_{N_s}(z) = {P^{(2N_s)}_1(z)}$ and 
$b_{N_s}(z) = z^{-8\ell_+}{P_2^{(2N_s)}(z)}$. The principal branch for the discrete 
scattering coefficients here works out to
be $\Re\zeta\in[-\pi/2h,\pi/2h]$. This again follows from the principle branch of the
individual transfer matrices. As before, the nodes
$\zeta_j=\xi_j+i0=j(\pi/h)/N,\,j\in\field{Z},$ lead to
$z_j=\exp(i\pi j/4N)$ which is not in the standard form for FFT algorithms to be
used. Therefore, we would like to work with $N'=8N$ nodes so that 
$z_j=\exp(i2\pi j/N')$ and pad the input vector with zeros.

\section{Numerical Tests and Conclusion}\label{sec:num-exp}
For the numerical experiments, we employ the well-known secant-hyperbolic 
potential given by 
$q(t) = A \sech{t}$, ($\kappa=-1$) for which the scattering coefficients are given
in~\cite{V2017INFT1}. We set the computational domain to be $[-30,30]$ and 
let $A\in\{4.4, 12.4\}$. Let $\Omega_h$ be the principal branch; then, the 
error in computing $b(\xi)$ is quantified by  
\begin{equation}
e_{\text{rel.}}=\|b(\xi)-b_N(\xi)\|_{\fs{L}^2(\Omega_h)}/\|b(\xi)\|_{\fs{L}^2(\Omega_h)},
\end{equation}
where the integrals are computed using the trapezoidal rule. Similar consideration applies 
to $\rho(\xi)$. For the purpose of testing, we include the implicit Adams method
presented in~\cite{V2018LPT} which are labelled as IA$_m$ with $m=2,3$. The
method IA$_1$ is identical to the trapezoidal rule, therefore, we use the label
TR. The convergence analysis is carried out in Fig.~\ref{fig:convg} and the
trade-off between accuracy and complexity is presented in Fig.~\ref{fig:trade-off}.
In terms of accuracy, the CF$_{24}$ outperforms every other method with M$_{34}$
being a close second. However, the accuracy-complexity trade-off is similar for the two
methods. The `fast' methods evidently lower complexity at the cost of accuracy.
The RK methods (ERK$_{34}$ and IRK$_{34}$) outperform all the other 'fast' methods 
in terms of accuracy-complexity trade-off (see Fig.~\ref{fig:trade-off}); however, with
increasing signal strength, the `slow' methods becoming equally competitive. In
fact, at moderately small step-sizes, the `slow' methods far outperform the
`fast' methods (see Fig.~\ref{fig:err}) with increasing signal strength.
\bibliographystyle{IEEEtran}
%\bibliography{IEEEabrv,NFT4_PTL}
\providecommand{\noopsort}[1]{}\providecommand{\singleletter}[1]{#1}%

\end{document}